\documentclass[12pt]{article}
\usepackage{amsfonts}
\usepackage[dvips]{graphics}
\newtheorem{theorem}{Theorem}

\newtheorem{definition}{Definition}
\newtheorem{remark}{Remark}
\newtheorem{example}{Example}

\makeindex
\begin{document}
\title{A  Liouville comparison principle  for weak solutions of semilinear parabolic
second-order partial differential inequalities in the whole space}
\author{Vasilii V. Kurta}
\maketitle
\thispagestyle{empty}
\begin{abstract}
\noindent \noindent We obtain a new Liouville comparison principle
for  weak solutions $(u,v)$ of semilinear parabolic second-order
partial differential inequalities of the form
$$
\qquad \qquad \qquad u_t -{\mathcal L}u- |u|^{q-1}u\geq v_t -{\mathcal L}v- |v|^{q-1}v\qquad \qquad \qquad   (*)
$$
in the whole space  ${\mathbb E} = {\mathbb R} \times \mathbb R^n$.
Here $n\geq 1$, $q>0$ and
$$
{\mathcal L}=\sum\limits_{i,j=1}^n\frac{\partial
}{{\partial}x_i}\left [ a_{ij}(t, x) \frac{\partial }{{\partial
}x_j}\right],$$ where   $a_{ij}(t,x)$, $i,j=1,\dots ,n$, are
functions that are defined, measurable and locally bounded in
$\mathbb E$, and such that $a_{ij}(t,x)=a_{ji}(t,x)$ and
$$ \sum_{i,j=1}^n a_{ij}(t,x)\xi_i\xi_j\geq 0
$$
for almost all $(t,x)\in \mathbb E$ and all $\xi \in \mathbb R^n$.
We show that the  critical exponents in the Liouville comparison
principle obtained, which are responsible for the non-existence of
non-trivial (i.e., such that $u\not \equiv v$)  weak solutions to
($*$) in the whole space $\mathbb E$, depend on the behavior of the
coefficients of the operator $\mathcal L$  at infinity and coincide
with those obtained for solutions of  ($*$) in the half-space
${\mathbb R}_+\times {\mathbb R}^n$. As direct corollaries we obtain
new Liouville-type theorems for non-negative weak solutions $u$ of
the inequality ($*$) in the whole space $\mathbb E$ in the case when
$v\equiv 0$. All the results obtained are new and sharp.
\end{abstract}

\section{Introduction and preliminaries}
This work should be considered as a supplement to the paper [5] and
is devoted to a new Liouville comparison principle for  weak
solutions to parabolic inequalities of the form
\begin{eqnarray}
 u_t -{\mathcal L}u - |u|^{q-1}u\geq v_t -{\mathcal L}v - |v|^{q-1}v
\end{eqnarray}
in the whole space  $\mathbb E= \mathbb R \times \mathbb R^n $,
where  $n\geq 1$ is a natural number,  $q>0$ is a real number and
${\mathcal L}$ is a linear second-order partial differential
operator  in divergence form defined   by the relation
\begin{eqnarray}
{\mathcal L}=\sum\limits_{i,j=1}^n\frac{\partial }{{\partial
}x_i}\left [ a_{ij}(t, x) \frac{\partial }{{\partial }x_j}\right]
\end{eqnarray}
for all $(t,x)\in \mathbb E$.  We assume that  the coefficients
$a_{ij}(t, x)$, $i,j=1, \dots, n$, of the operator ${\mathcal L}$
are functions that are defined, measurable and locally bounded  in
$\mathbb E$. We also assume that $a_{ij}(t, x)=a_{ji}(t, x)$,
$i,j=1,\dots ,n$, for almost all $(t, x)\in \mathbb E$, and  that
the corresponding quadratic form satisfies the conditions
\begin{equation} 0\leq \sum\limits_{i,j=1}^n
a_{ij}(t,x)\xi_i\xi_j\leq A(t, x)|\xi|^2 \end{equation} for all
$\xi=(\xi_1, \dots, \xi_n)\in {\mathbb R}^n$ and  almost all $(t,
x)\in \mathbb E$, with $A(t, x)$ a  function that is defined,
measurable, non-negative  and locally bounded in $\mathbb E$.

It is important to note that if $u=u(t,x)$  satisfies the inequality
\begin{eqnarray} u_t \geq {\mathcal L}u+|u|^{q-1}u
\end{eqnarray} and  $v=v(t,x)$ satisfies the inequality
\begin{eqnarray} v_t \leq {\mathcal
L}v+|v|^{q-1}v,\end{eqnarray}  then the pair $(u,v)$ satisfies the
inequality (1). Thus, all  the results obtained in this paper for
solutions of (1) are valid for the corresponding solutions of the
equations (4), (5).

The  results obtained in [5] for solutions of the inequality (1) in
the half-space $\mathbb S= (0, +\infty) \times \mathbb R^n$, $n\geq
1$,  show that the behavior of the coefficients $a_{ij}(t,x)$ of the
operator $\mathcal L$ as $|x|\to + \infty$  manifests itself in
Liouville-type results; namely, the  critical exponents in the
Liouville comparison principle for  weak solutions of (1) in the
half-space  $\mathbb S$, which are responsible for the non-existence
of non-trivial (i.e., such that $u\not\equiv v$)  weak solutions to
the inequality (1) in  $\mathbb S$, depend in an essential way on
the behavior of the coefficients of the operator $\mathcal L$ as
$|x|\to + \infty$.

The main goal of the present  work is to show  that  similar
critical exponents in a Liouville comparison principle for  weak
solutions of (1) in the whole space $\mathbb E$ also exist and, what
is  more intriguing,  coincide with those obtained in [5] for
solutions of  (1) in the half-space $\mathbb S$. In this connection
it is important to note here that the latter, generally speaking, is
not the case for solutions of equations corresponding to the
inequalities (1), (4) and (5). To make certain of this, it is enough
to compare the famous Fujita critical blow-up exponent $q_F=1+2/n$
for non-negative classical solutions of the equation
\begin{eqnarray} u_t ={\Delta }u+|u|^{q-1}u
\end{eqnarray}
in the half-space $\mathbb S$ obtained in [2--4] with the blow-up
exponent for non-negative classical solutions of the equation (6) in
the whole space $\mathbb E$ obtained in [1] and [8], which is equal
to
\begin{eqnarray}
q_{B}=\left \{ \begin {array} {ll}\frac{n(n+2)}{(n-1)^2},
&
{\rm if}\   n\geq 2,\\
+\infty, &
{\rm if}\ n=1.\end{array}\right.\nonumber
\end{eqnarray}

In order to trace the relation between the behavior of the
coefficients $a_{ij}(t,x)$ of the operator $\mathcal L$ as $|x|\to +
\infty$ and the critical exponents which are responsible for the
non-existence of non-trivial   weak solutions to the inequality (1)
in the whole space $\mathbb E$, we consider the quantity
\begin{eqnarray}{\mathcal A}(R)={\mathrm {ess}\sup}_{(t,x)\in
(-\infty,+\infty )\times \{R/2 <|x|<R \}}A(t,x)\end{eqnarray} for any
$R>0$ and  assume that the coefficients of the operator $\mathcal L$
satisfy  the condition
\begin{eqnarray}{\mathcal A}(R)\leq cR^{2-\alpha},\end{eqnarray}
with some real constant $\alpha$ and some real positive constant
$c$, for all $R>1$. It is clear that if $\alpha <2$, then the
coefficients of the operator $\mathcal L$ may be unbounded in
$\mathbb E$, if $\alpha =2$, they are globally bounded in $\mathbb
E$, and if $\alpha
>2$, they must vanish as $|x|\to +\infty$.

We also introduce a special function space, which is directly
associated with the linear partial differential operator
$\displaystyle {\mathcal P}= \frac{\partial}{\partial t}- {\mathcal
L}, $ and assume that the weak solutions of the inequalities (1),
(4) and (5) belong to this function space only locally  in $\mathbb
E$.

\section{Definitions}

Let ${\mathbb S}_\tau=(\tau, +\infty)\times {\mathbb R}^n$, where
$\tau$ is a real number or $\tau=-\infty$, and $n\geq 1$. It is
clear  that the half-space $\mathbb S$ and the whole space $\mathbb
E$ are partial cases of the set ${\mathbb S}_\tau$ when $\tau=0$ and
$\tau=-\infty$.

\begin{definition}
Let   $n\geq 1$, $q>0$,  $\hat q =\max\{1,q\}$ and $\tau\in
[-\infty, +\infty)$, let $\mathcal L$ be a differential operator
defined by (2) in the whole space $\mathbb E$, and let $\Omega$ be
an arbitrary domain in ${\mathbb S}_\tau$. By $W^{\mathcal L,
q}(\Omega)$ we denote the completion of the function space
$C^{\infty}(\Omega)$ with respect to the norm
\begin{eqnarray}
\| w \|_{W^{\mathcal L, q}(\Omega)}= \int\limits_{\Omega} |w_t|dtdx+
\left[\int\limits_{\Omega}\sum\limits_{i,j=1}^n
a_{ij}(t,x)\frac{\partial w}{\partial x_i}\frac{\partial w}{\partial
x_j}dtdx\right]^{1/2} + \left[\int\limits_{\Omega}|w|^{\hat
q}dtdx\right]^{1/{\hat q}}\nonumber
\end{eqnarray}
where $C^{\infty}(\Omega)$ is the space of all functions  defined
and infinitely differentiable  in  $\Omega$.
\end{definition}
\begin{definition}
Let  $n\geq 1$,  $q>0$ and $\tau\in [-\infty, +\infty)$, and let
$\mathcal L$ be a differential operator defined  by (2) in the whole
space $\mathbb E$. A function $w=w(t,x)$ belongs to the function
space $W^{\mathcal L, q}_{\mathrm {loc}}({\mathbb S}_\tau)$ if $w$
belongs to $W^{\mathcal L, q}(\Omega)$ for any bounded domain
$\Omega$ in ${\mathbb S}_\tau$.
\end{definition}
\begin{definition}
Let  $n\geq 1$, $q>0$ and $\tau\in [-\infty, +\infty)$, and let
$\mathcal L$ be a differential operator defined  by (2) in the whole
space $\mathbb E$.  A pair $(u,v)$ of functions $u=u(t,x)$ and
$v=v(t,x)$ is called  a  weak solution to the inequality (1) in
${\mathbb S}_\tau$, if these functions are defined and measurable in
${\mathbb S}_\tau$, belong to the function space $W^{\mathcal L,
q}_{\mathrm {loc}}({\mathbb S}_\tau)$ and satisfy the integral
inequality
\begin{eqnarray}
\int\limits_{{\mathbb S}_\tau}\left[u_t\varphi+\sum_{i,j=1}^n
a_{ij}(t,x)\frac{\partial \varphi}{\partial x_i} \frac{\partial
u}{\partial x_j} - |u|^{q-1}u\varphi\right]dtdx \geq \nonumber
\\
\int\limits_{{\mathbb S}_\tau}\left[v_t\varphi+\sum_{i,j=1}^n
a_{ij}(t,x)\frac{\partial \varphi}{\partial x_i} \frac{\partial
v}{\partial x_j}- |v|^{q-1}v\varphi\right]dtdx
\end{eqnarray} for
every  function $\varphi \in  C^\infty ({\mathbb S}_\tau)$ with
compact support in ${\mathbb S}_\tau$,  where $C^{\infty}({\mathbb
S}_\tau)$ is the space of all functions defined and infinitely
differentiable in ${\mathbb S}_\tau$.
\end{definition}

\begin{remark} We understand the inequality  (9) in the sense
discussed, e.g., in [7].
\end{remark}

Analogous definitions of solutions to the inequality  (4) and  the
inequality  (5) in ${\mathbb S}_\tau$, as special cases of the
inequality  (1) in ${\mathbb S}_\tau$ for $v\equiv 0$ or $u\equiv
0$, follow immediately from Definition 3.

It is clear that if $\tau_1<\tau_2$, then the  solutions of the
inequality (1), (4) or (5) in ${\mathbb S}_{\tau_1}$ are solutions
of the corresponding inequality  in ${\mathbb S}_{\tau_2}$. In
particular, if $(u,v)$ is a solution of the inequality (1) in
$\mathbb E$, then $(u,v)$ is  a solution of the inequality (1) in
${\mathbb S}_{\tau}$ for any $\tau\in (-\infty, +\infty)$.

\section{Results}
%Theorem 1
\begin{theorem} Let $n\geq 1$, $\alpha>0$ and  $1<q\leq 1 +\frac \alpha
n$, let $\mathcal L$ be a differential operator defined  by (2) in
the whole space $\mathbb E$, the coefficients of which  satisfy the
condition (8) with the given $\alpha$ and some $c>0$,   and  let
$(u,v)$ be a  weak solution of the inequality (1) in  $\mathbb E$
such that $u\geq v$. Then $u = v$ in $\mathbb E$.
\end{theorem}

As we have observed above, since any pair of solutions $u=u(t,x)$,
$v=v(t,x)$ of the inequalities (4), (5) in $\mathbb E$ is  a
solution $(u,v)$  of the inequality (1) in  $\mathbb E$,  the
following statement is a direct corollary of Theorem 1.
%Theorem 2
\begin{theorem} Let $n\geq 1$, $\alpha>0$ and  $1<q\leq 1 +\frac \alpha
n$, let $\mathcal L$ be a differential operator defined  by (2) in
the whole space $\mathbb E$, the coefficients of which   satisfy the
condition (8) with the given $\alpha$ and some $c>0$,   and let
$u=u(t,x)$ be a  weak solution of the inequality (4) and  $v=v(t,x)$
be a  weak solution of the inequality (5) in  $\mathbb E$ such that
$u\geq v$. Then $u = v$ in $\mathbb E$.
\end{theorem}

Each of the results in Theorems 1 and 2, which obviously  have the
character of a comparison principle, we term a Liouville-type
comparison principle, since in  particular cases when either
$u\equiv 0$ or $v\equiv 0$, it becomes a Liouville-type theorem for
solutions of (5) or (4), respectively. We formulate here only  the
case when $v\equiv 0$.

%Theorem 3
\begin{theorem} Let $n\geq 1$, $\alpha>0$ and  $1<q\leq 1 +\frac \alpha
n$, let $\mathcal L$ be a differential operator defined  by (2)  in
the whole space $\mathbb E$, the coefficients of which   satisfy the
condition (8) with the given $\alpha$ and some $c>0$,   and let
$u=u(t,x)$ be a non-negative
 weak solution of the inequality (4) in  $\mathbb E$. Then $u =
0$ in $\mathbb E$.
\end{theorem}

Note that all the results in Theorems 1--3, including the partial
case when  $\mathcal L$ is  the Laplacian operator, are new and
sharp. (We demonstrate their sharpness below by Examples 1--3).
Thus, as we have already mentioned above, the critical exponents in
Theorems 1--3,  which are  responsible for the non-existence of
non-trivial  weak solutions to the inequalities (1), (4) and (5) in
the whole space $\mathbb E$, coincide with those obtained in
Theorems 1--3 in [5] for weak solutions of the corresponding
inequalities in the half-space $\mathbb S$. In the particular case
when $\alpha=2$, the critical exponent in Theorems 1--3  coincides
with the well-known Fujita critical blow-up exponent obtained in
[2--4].

%Example 1
\begin{example}
Let $n\geq 1$, $ \alpha>-\infty $ and $q\leq 1$, and let $\mathcal
L$ be a differential operator defined  by (2) in the whole space
$\mathbb E$, the coefficients of which  satisfy the condition (8)
with the given $\alpha$ and some $c>0$. It is clear that the
function $u(t,x)=\exp(t)$ is a positive classical solution of the
inequality (4) in the whole space $\mathbb E$. Also, it is clear
that  the function $v=-u(t,x)$ is a negative classical solution of
the inequality (5) in the whole space $\mathbb E$. Thus, the pair of
functions $u=u(t,x)$ and $v=v(t,x)$ is a non-trivial classical
solution of the system (4)--(5) and, therefore, $(u,v)$ is a
non-trivial classical solution of the inequality (1) in the whole
space $\mathbb E$ such that $u(t,x)> v(t,x)$.
\end{example}

%Example 2
\begin{example} Let  $n\geq 1$,  $\alpha>0$ and  $q> 1 +\frac \alpha n$.
Consider the operator $\mathcal L$ defined by (2) in the whole space
$\mathbb E$ with the coefficients given by the expression
\begin{eqnarray}
a_{ij}(t,x)=(1+|x|^2)^{\frac{2-\alpha}2}\delta_{ij}
\end{eqnarray}
for all $(t,x)\in \mathbb E$, where $\delta_{ij}$ are Kronecker's
symbols and $i,j=1,\dots ,n$. It is easy to see that the condition
(8) is fulfilled for these coefficients with the given $\alpha$ and
some $c>0$. Further, for the given $\alpha$,  let
\begin{eqnarray}
u(t,x)=\left \{ \begin {array} {ll}\kappa t^{-\beta}{\mathcal E}(t,x),
&
{\rm if}\   t>0,\, x\in {\mathbb R}^n,\\
0, &
{\rm if}\   t\leq 0,\, x\in {\mathbb R}^n,\end{array}\right.
\end{eqnarray}
where ${\mathcal E}(t,x)=\exp\left(-\gamma \frac{(1+|x|^2)^{\frac
\alpha 2}}{t}\right)$ for all $t> 0$ and $x\in {\mathbb R}^n$ and
the positive constants  $\beta$, $\gamma$ and $\kappa$ will be
chosen  below.

First, since the function $u=u(t,x)$ of the form (11) with any fixed
positive constants  $\alpha$, $\beta$, $\gamma$ and $\kappa$ is
infinitely differentiable in the whole space $\mathbb E$ and
vanishes, along with  all its derivatives, for all $t\leq 0$ and
$x\in {\mathbb R}^n$, it is clear that  $u=u(t,x)$ is a classical
solution of the inequality (4)  for all $t\leq 0$ and $x\in {\mathbb
R}^n$.

Now, consider the case when $t>0$ and $x\in {\mathbb R}^n$. Making
necessary calculations, we have
\begin{eqnarray}
u_t= -\kappa\beta t^{-\beta -1}{\mathcal E}(t,x)+\kappa \gamma t^{-\beta -2}(1+|x|^2)^{\frac{\alpha}2}{\mathcal E}(t,x),\nonumber
\end{eqnarray}
\begin{eqnarray}
\frac{\partial u}{\partial x_i}=-\alpha \kappa\gamma  t^{-\beta
-1}(1+|x|^2)^{\frac{\alpha}2-1}{\mathcal E}(t,x)x_i\nonumber
\end{eqnarray} and
\begin{eqnarray}
\frac{\partial }{\partial x_i}\left(a_{ii}(t,x)\frac{\partial u}{\partial x_i}\right)=\nonumber \\
-\alpha \kappa\gamma  t^{-\beta -1}{\mathcal E}(t,x)+ \alpha^2 \kappa \gamma^2
t^{-\beta -2}(1+|x|^2)^{\frac{\alpha}2}\frac{{x_i}^2}{1+|x|^2}{\mathcal E}(t,x),\nonumber
\end{eqnarray}
for all $t> 0$ and $x\in {\mathbb R}^n$, where  the coefficients
$a_{ii}(t,x)$ are given by (10) and $i=1,\dots, n$. Further, it is
also easy to calculate that
\begin{eqnarray}
u_t-{\mathcal L}u=\nonumber \\
(\alpha \kappa n \gamma -\kappa \beta)  t^{-\beta -1}{\mathcal E}(t,x) + \left(\kappa \gamma - \alpha^2 \kappa \gamma^2\frac{|x|^2}{1+|x|^2}\right)
t^{-\beta -2}(1+|x|^2)^{\frac{\alpha}2}{\mathcal E}(t,x)\nonumber
\end{eqnarray} and
\begin{eqnarray}
|u|^{q-1}u= {\kappa}^q t^{-\beta q }{\mathcal E}^q(t,x).\nonumber
\end{eqnarray}

As a result, the inequality (4)  with $u=u(t,x)$ given by (11) takes
the form
\begin{eqnarray}
(\alpha \kappa n \gamma -\kappa \beta)  t^{-\beta -1}{\mathcal E}(t,x) +\nonumber\\ \left(\kappa \gamma - \alpha^2 \kappa \gamma^2\frac{|x|^2}{1+|x|^2}\right)
t^{-\beta -2}(1+|x|^2)^{\frac{\alpha}2}{\mathcal E}(t,x)\geq {\kappa}^q t^{-\beta q }{\mathcal E}^q(t,x)
\end{eqnarray}
for all $t>0$ and $x\in {\mathbb R}^n$.  Now, choosing the constants
$\beta$, $\gamma$ and $\kappa$ such that
\begin{eqnarray}
\begin {array} {ll} \beta=\frac 1{q-1}, \quad \frac 1{\alpha n(q-1)}<\gamma\leq \left(
\frac 1{\alpha}\right)^2, \\ 0< \kappa\leq\left ( \alpha
n\left(\gamma -\frac 1{\alpha
n(q-1)}\right)\right)^{1/(q-1)}\end{array}\end{eqnarray} and taking
into account that ${\mathcal E}(t,x)\leq 1$ for all $t> 0$ and $x\in
{\mathbb R}^n$, it is not difficult to verify that the inequality
(12) holds  for all $t> 0$ and $x\in {\mathbb R}^n$, and, therefore,
the function $u=u(t,x)$ of the form (11), with the given $\alpha$
and $q$, and the constants $\beta$, $\gamma$ and $\kappa$ satisfying
the conditions (13), is  a positive classical solution of the
inequality (4) for all $t> 0$ and $x\in {\mathbb R}^n$.

Thus, we may conclude  that the function $u=u(t,x)$ of the form
(11), with the given $\alpha$ and $q$, and the constants $\beta$,
$\gamma$ and $\kappa$ satisfying the conditions (13), is indeed a
non-trivial non-negative classical solution of the inequality (4) in
the whole space $\mathbb E$, with $a_{ij}(t,x)$, the coefficients of
the operator $\mathcal L$, defined  by (10). Also, it is clear that
the function $v=-u(t,x)$ is a non-trivial non-positive classical
solution of the inequality (5) in the whole space $\mathbb E$, with
$a_{ij}(t,x)$ in (2) defined by (10). Thus, the pair  of  functions
$u=u(t,x)$ and $v=v(t,x)$ is a non-trivial classical solution of the
system (4)--(5) and, therefore, $(u,v)$ is a non-trivial classical
solution of the inequality (1) in the whole space $\mathbb E$ such
that $u(t,x)\geq v(t,x)$, with $a_{ij}(t,x)$ in (2) defined by (10).
\end{example}

Note that non-negative  classical super-solutions to linear
uniformly pa\-ra\-bolic equations with globally bounded coefficients
in non-divergence form in the whole space $\mathbb E$ except the
origin of coordinates  in a form close to that given by the relation
(11) with $\alpha=2$ were  constructed in [6, p. 122]. Also, note
that positive classical super-solutions of the equation (6) in the
half-space $\mathbb S$ in a form close to that given by the relation
(11) with $\alpha=2$ were constructed in [9, p. 283].

%Example 3
\begin{example} Let  $n\geq 1$, $\alpha\leq 0$,  $q>1 +\frac \alpha
n$ and $q>1$,   and let $\hat \alpha$ be any positive number such
that $q>1+\frac {\hat \alpha}n$. Consider the operator $\mathcal L$
defined by (2) in the whole space  $\mathbb E$ with the coefficients
given by the relation
\begin{eqnarray}
a_{ij}(t,x)=(1+|x|^2)^{\frac{2-\hat \alpha}2}\delta_{ij}
\end{eqnarray}
for all $(t,x)\in \mathbb E$, where $\delta_{ij}$ are Kronecker's
symbols and $i,j=1,\dots,n$. As in Example 2, it is easy to see that
$\displaystyle{\mathcal A}(R)\leq CR^{2-\hat\alpha}$ for all $R>1$,
with $C$ some positive constant which possibly depends on $\hat
\alpha$ and $n$, and, therefore,  the condition (8) is fulfilled for
these coefficients with the given $\alpha$ and some $c>0$. Also, for
the given $\hat \alpha$ and $q$, let $\beta=\frac 1{q-1}$, $\frac
1{\hat \alpha n(q-1)}<\gamma\leq \left (\frac 1{\hat
\alpha}\right)^2$, $0< \kappa\leq\left ( \hat \alpha n\left(\gamma
-\frac 1{\hat \alpha n(q-1)}\right)\right)^{1/(q-1)}$ and
\begin{eqnarray}
u(t,x)=\left \{ \begin {array} {ll}\kappa t^{-\beta}\exp\left(-\gamma \frac{(1+|x|^2)^{\frac
{\hat\alpha} 2}}{t}\right),
&
{\rm if}\   t>0,\, x\in {\mathbb R}^n,\\
0, &
{\rm if}\   t\leq 0,\, x\in {\mathbb R}^n.\end{array}\right.
\end{eqnarray}

Again, as in Example 2, it is not difficult to verify that the
function $u=u(t,x)$ defined by the expression (15) is a nontrivial
non-negative classical  solution of the inequality (4) in the whole
space $\mathbb E$, with $a_{ij}(t,x)$ in (2) defined by (14). Also,
it is clear that the function $v=-u(t,x)$ is a nontrivial
non-positive classical solution of the inequality (5) in the whole
space $\mathbb E$, with $a_{ij}(t,x)$, the coefficients of the
operator $\mathcal L$, defined by (14). Thus, the pair of functions
$u=u(t,x)$ and $v=v(t,x)$ is a non-trivial classical solution of the
system (4)--(5) and, therefore, $(u,v)$ is a non-trivial classical
solution of the inequality (1) in the whole space $\mathbb E$ such
that $u(t,x)\geq  v(t,x)$, with $a_{ij}(t,x)$ in (2) defined by
(14).
\end{example}

%Remark
\begin{remark} For the case when $n\geq 1$,  $\alpha\leq 0$ and
$1 \geq q> 1 +\frac \alpha
n$, see Example 1.
\end{remark}

Note that Examples 1--3 are constructed on the basis of those in
[5].

\section{ Proofs}

First of all, we would like to underline that one can prove Theorem
1 by following the approach  in [5] with easy modifications. Here,
however, to  prove  Theorem 1 we use  different arguments which are
based on the structure of the conditions (7) and (8) on the behavior
of the coefficients of the operator $\mathcal L$ as $|x|\to
+\infty$, in particular, on the invariance of the quantity $\mathcal
A(R)$ in (7) with respect to the shifts along the $t$-axis, and
technically which consist of  applying the Liouville comparison
principle for weak solutions of the inequality (1) in the half-space
$\mathbb S$ obtained in [5]. For the benefit of the reader we now
formulate  that result from [5] in the following form:

{\bf Theorem 1$'$} {\it Let $n\geq 1$, $\alpha>0$ and  $1<q\leq 1
+\frac \alpha n$, let  ${\mathcal L}$ be a differential operator
defined by (2) in the whole space $\mathbb E$, the coefficients of
which are equal to zero on the set $\mathbb E\setminus \mathbb S$
and satisfy the condition (8) with the given $\alpha$ and some
$c>0$, and  let $(u,v)$ be a weak solution of the inequality  (1) in
the half-space $\mathbb S$ such that $u\geq v$. Then $u = v$ in
$\mathbb S$.}

{\it Proof of Theorem 1.} The proof is by contradiction. Let $n\geq
1$, $\alpha>0$ and  $1<q\leq 1 +\frac \alpha n$, let $\mathcal L$ be
a differential operator defined  by (2) in the whole space $\mathbb
E$ with  the coefficients  $a_{ij}(t,x)$ satisfying  the condition
(8) with the given $\alpha$ and some $c>0$, and  let $(u,v)$ be a
nontrivial  weak solution of the inequality (1) in $\mathbb E$ such
that $u\geq v$. Then there exists a real number $\tau$ such that
$u\not \equiv v$ in ${\mathbb S}_\tau$. (If such a number $\tau$
does not exist, then the solution $(u,v)$ of (1) in $\mathbb E$ is
trivial, i.e., $u(t,x) = v(t,x)$ almost everywhere in $\mathbb E$.)
Therefore, the pair $(u,v)$ is, in particular, a non-trivial
solution of the inequality (1) in ${\mathbb S}_\tau$ such that $
u\geq v$.

Further, consider the pair $(\hat u, \hat v)$ of functions $\hat
u=\hat u(t,x)$ and $ \hat v=\hat v(t,x)$ given by the relations
\begin{eqnarray}
\hat u(t,x)=u (t+\tau, x), \quad \hat v(t,x)=v  (t+\tau, x)
\end{eqnarray}
for all $(t,x)\in \mathbb E$. It is clear that the pair $(\hat
u,\hat v)$ is a non-trivial, i.e.,
\begin{eqnarray}
 \hat u\not \equiv \hat v,
\end{eqnarray}
solution of the inequality
\begin{eqnarray}
 {\hat u}_t -\hat{\mathcal L}\hat u - |\hat u|^{q-1}\hat u\geq {\hat v}_t -\hat {\mathcal L}\hat v - |\hat v|^{q-1}\hat v\nonumber
\end{eqnarray}
of the form (1) in  $\mathbb S$ such that $\hat u\geq \hat v$, with
the differential operator
\begin{eqnarray}
\hat {\mathcal L}=\sum\limits_{i,j=1}^n\frac{\partial
}{{\partial}x_i}\left [ \hat{a}_{ij}(t, x) \frac{\partial
}{{\partial }x_j}\right] \nonumber \end{eqnarray}
of the form (2), where
\begin{eqnarray}
\hat{a}_{ij}(t,x)= \left \{ \begin {array} {ll} {a}_{ij}(t+\tau,x),
&
{\rm if} \   t>0,\, x\in {\mathbb R}^n,\\
0, & {\rm if} \  t\leq 0,\, x\in {\mathbb R}^n,\end{array}\right.
\nonumber \end{eqnarray} with  $i,j=1, \dots, n$. It is also clear
that the inequalities
\begin{eqnarray} 0\leq \sum\limits_{i,j=1}^n
\hat{a}_{ij}(t,x)\xi_i\xi_j\leq A(t+\tau, x)|\xi|^2 \nonumber
\end{eqnarray} hold  for all $\xi=(\xi_1, \dots, \xi_n)\in
{\mathbb R}^n$ and almost all $(t, x)\in \mathbb E$, where  $A(t,
x)$ is the function from the conditions (3).

Now, consider the quantity
\begin{eqnarray}\hat {\mathcal A}(R)={\mathrm {ess}\sup}_{(t,x)\in
(-\infty,+\infty )\times \{R/2 <|x|<R \}}A(t+\tau,x)\end{eqnarray}
for any $R>0$. It follows easily from (7) and (18) that
\begin{eqnarray}\hat {\mathcal A}(R)=  {\mathcal A}(R)\nonumber \end{eqnarray}
for any $R>0$, and, therefore, the coefficients of the operator
$\hat{\mathcal L}$ and the coefficients of the operator ${\mathcal
L}$ satisfy the condition (8) with exactly the same constants $c$
and $\alpha$.

As a result, for $n\geq 1$, $\alpha>0$ and  $1<q\leq 1 +\frac \alpha
n$, we have the pair $(\hat u, \hat v)$ consisting of the functions
$\hat u=\hat u(t,x)$ and $ \hat v=\hat v(t,x)$ given by the
relations (16), which is a weak solution of the inequality (1) in
the half-space $\mathbb S$ such that $\hat u\geq \hat v$, with
${\mathcal L}=\hat {\mathcal L}$. From our consideration, it is easy
to see that the coefficients of the operator $\hat {\mathcal L}$
satisfy all the hypotheses of Theorem 1$'$ and that the constants
$\alpha$ and $c$ in the formulations of Theorem 1 and Theorem 1$'$
coincide. Then, by Theorem 1$'$, we have to conclude that $\hat
u(t,x)=\hat v(t,x)$ almost everywhere in the half-space $\mathbb S$,
which  in turn contradicts the relation (17). Thus, the
contradiction obtained proves Theorem 1.

\vspace{10mm}

\newpage

\noindent \textbf{Authors' addresses:}

\vspace{5 mm}

\noindent Vasilii V. Kurta

\noindent Mathematical Reviews

\noindent 416 Fourth Street, P.O. Box 8604

\noindent Ann Arbor, Michigan 48107-8604, USA

\noindent \textbf {e-mail:} vkurta@umich.edu, vvk@ams.org
\end{document}